\documentclass[9pt,reqno]{amsart}
\usepackage{graphicx}
\usepackage{verbatim}
\usepackage{textcomp}
\usepackage{amssymb}
\usepackage{cite}
\usepackage{amsmath}
\usepackage{latexsym}
\usepackage{amscd}
\usepackage{amsthm}
\usepackage{mathrsfs}
\usepackage{bm}
\usepackage{url}
\usepackage{hyperref}
\usepackage{bookmark}
\allowdisplaybreaks[3]
\newtheorem{thm}{Theorem}[section]

\newtheorem{lem}[thm]{Lemma}

\theoremstyle{definition}

\theoremstyle{remark}

\numberwithin{equation}{section}
\begin{document}
\title[An overdetermined problem]
{An overdetermined problem related to the $p$-Laplacian on Riemannian manifolds}

\author{Guangyue Huang}
\author{Chunlei Luo}
\author{Hongru Song$^*$}

\address{School of Mathematics and Statistics, Henan Normal University, Xinxiang 453007, P.R.China  }

\email{hgy@htu.edu.cn(G. Huang) }
\email{clluo2024@126.com(C. Luo) }

\email{songhongru@htu.edu.cn(H. Song) }

\thanks{The research of the first author is supported by key projects of the Natural Science Foundation of Henan Province (No. 252300421303) and NSFC (No. 11971153), the third author is supported by Postdoctoral Fellowship Program of CPSF (No. GZC20230734).}
\thanks{$^*$Corresponding author: songhongru@htu.edu.cn(H. Song) }

\begin{abstract}
In this paper, we study the overdetermined problem for the $p$-Laplacian equation on a compact Riemannian manifold
with positive Ricci curvature. By introducing a new $P$-function which is related to the first nonzero eigenvalue for  $p$-Laplacian, we obtain some integral identities. As their applications, the Heintze-Karcher type inequality and the Soap Bubble Theorem have been achieved.

\end{abstract}

\subjclass[2020]{35A23; 53C24.}

\keywords{$p$-Laplacian; overdetermined problem; rigidity}

\maketitle

\section{Introduction}

Let $\Omega\subset \mathbb{R}^n$ be a bounded domain with boundary $\partial \Omega$, Serrin \cite{Serrin1971} studied the rigidity of the following overdetermined boundary value problem:
\begin{equation}\label{1-1}
\Delta u=-1,\ {\rm in}\ \Omega,
\end{equation}
\begin{equation}\label{1-2}
u=0,\ \ \ u_{\nu}=c,\ \   {\rm on}\ \partial \Omega,
\end{equation}
where $c$ is a constant and $\nu$ denotes the outward unit normal on $\partial \Omega$. Employing the method of moving planes, Serrin proved that the existence of a solution to \eqref{1-1}-\eqref{1-2} forces $\Omega$ to be a ball and $u$ to be radial. Subsequently, Weinberger gave a simple proof using integral identities, now known as the $P$-function approach. Specifically, Weinberger introduced a sub-harmonic function $P$. An application of the classical strong maximum principle then shows that $P$ is constant, implying the rigidity of the problem. Unlike Serrin's approach, Weinberger's argument crucially depends on the linearity of the Laplacian operator. Following the seminal works of Serrin and Weinberger, numerous generalizations of such overdetermined problems have been established for general elliptic operator in $\mathbb{R}^n$. For a comprehensive overview, we refer the interested reader to \cite{Birindelli2013,Buttazzo2011,Cianchi2009,Farina2010,Garofalo1989,WX2011} and the references therein.

It is also natural to consider a generalization of Serrin's symmetry result to quasilinear elliptic operators. The well-known $p$-Laplacian is defined by
\begin{equation}\label{1-3}
\Delta_p u={\rm div}(|\nabla u|^{p-2}\nabla u)
\end{equation}
with $p>1$. Let $\Omega\subset \mathbb{R}^n$ be a bounded domain. If the following overdetermined problem
\begin{equation}\label{1-4}
\Delta_p u=-1,\ {\rm in}\ \Omega,
\end{equation}
\begin{equation}\label{1-5}
u=0,\ \ \ u_{\nu}=c,\ \   {\rm on}\ \partial \Omega,
\end{equation}
admits a weak solution, then $\Omega$ must be a ball.
The result for $p$-Laplacian overdetermined problem \eqref{1-4}-\eqref{1-5} was first proved by Garofalo and Lewis \cite{Garofalo1989} via Weinberger's approach. For the case of $1<p<2$, Damascelli and Pacella \cite{Damascelli2000} proved this result by virtue of Serrin's approach. When $p>2$, Brock and Henrot \cite{Brock2002} gave a different proof of the result by using Steiner symmetrization. For the Serrin-type overdetermined problems involving the $k$-Hessian equation, see \cite{Brandolini2008,Ciraolo2017,Ciraolo2019,DY2022,Farina2022,GMY2023,GJZ2023,Ruan2024} and
the references therein.

As far as we know, there are few results on Riemannian manifolds with positive Ricci curvature by using the $P$-function approach since the $p$-Laplacian is nonlinear. In this paper, we introduce the so-called $P$-function associated with the first nonzero eigenvalue of the $p$-Laplacian and prove that this function is sub-harmonic with respect to the linearized operator of the $p$-Laplacian.

Let $(M,g)$ be an $n$-dimensional compact Riemannian manifold with boundary $\partial M$ and the Ricci curvature satisfying ${\rm Ric}\geq(n-1)K$, where $K$ is a positive constant. The mean curvature $H$ of $\partial M$ is given by $H=\frac{1}{n-1}{\rm tr}_g(II)$, where $II(X, Y)=g(\nabla_X\nu, Y)$ denotes the second fundamental form of $\partial M$ with $\nu$ the outward unit normal on $\partial M$. Let $u$ be a weak solution to
\begin{equation}\label{1-6}
\Delta_{p} u+\lambda u^{p-1}=-1,\ u>0 \ \ {\rm in}\ M,
\end{equation}
\begin{equation}\label{1-7}
u=0,\ \   {\rm on}\  \partial M,
\end{equation}
with the parameter $\lambda$ chosen as
$$
\lambda=\dfrac{(nK)^{\frac{p}{2}}}{(p-1)^{p-1}}
$$
because $\lambda$ is the lower bound of the first nonzero eigenvalue for  $p$-Laplacian\cite{WL2016}.
Therefore, we introduce the so-called $P$-function corresponding to \eqref{1-6} as follows:
\begin{align}\label{1-8}
P=(p-1)|\nabla u|^{p}+bu^p+\dfrac{p}{n}u,
\end{align}
where
$$
b=\Big[\dfrac{1}{n}-\Big(1-\dfrac{1}{n}\Big)\Big(1-\dfrac{2}{p}\Big)\Big]\lambda.
$$
For a positive weak solution $u$ to \eqref{1-6}, the outward normal $\nu$ is given by
\begin{equation}\label{1-9}
\nu=-\dfrac{\nabla u}{|\nabla u|}.
\end{equation}

Firstly, we obtain the following integral identities:

\begin{thm}\label{1-th-1}
Let $(M,g)$ be an $n$-dimensional compact Riemannian manifold with boundary $\partial M$ and the Ricci curvature satisfying ${\rm Ric}\geq(n-1)K$, where $K$ is a positive constant. Suppose $u$ is a weak solution to \eqref{1-6}, and let $\mathcal{L}_{u}$ be the operator defined in \eqref{2-1}. Then for $p\geq2$, we have
\begin{align}\label{1-th-1-F1}
\dfrac{n}{p(p-1)(n-1)}\int_{M} \mathcal{L}_{u}Pdv=|M|+\lambda\int_{M} u^{p-1}dv-\int_{\partial M} nH|u_{\nu}|^{2p-2}ds,
\end{align}
where $|M|$ denotes the volume of $M$.
\end{thm}

\begin{thm}\label{1-th-2}
Let $(M,g)$ be an $n$-dimensional compact Riemannian manifold with the boundary $\partial M$ and the Ricci curvature satisfying ${\rm Ric}\geq(n-1)K$, where $K$ is a positive constant. If the mean curvature $H$ is positive, then for the weak solution $u$ to \eqref{1-6} with $p\geq2$, we have
\begin{align}\label{1-th-2-F1}
\dfrac{n}{p(p-1)(n-1)}&\int_{M}\mathcal{L}_{u}Pdv+\int_{\partial M} \dfrac{1}{nH}(1+nH|u_{\nu}|^{p-2}u_{\nu})^2ds\notag\\
=&\int_{\partial M} \dfrac{1}{nH}ds-\lambda\int_{M} u^{p-1}dv-|M|.
\end{align}
\end{thm}

In the case of manifolds with nonnegative Ricci curvature, we have $K=0$ and consequently $\lambda=0$. Then the $P$-function defined in \eqref{1-8} reduces to $P=(p-1)|\nabla u|^{p}+\dfrac{p}{n}u$, which coincides with formula (5) in \cite{Ruan2024}. For this setting, the analogues of Theorems \ref{1-th-1} and \ref{1-th-2} have been proved by Ruan, Huang and Chen (see Theorem 2 and Theorem 3 in \cite{Ruan2024}).

Formula \eqref{1-th-2-F1} implies the following Heintze-Karcher type inequality:

\begin{thm}\label{1-th-3}
Let $(M,g)$ be an $n$-dimensional compact Riemannian manifold with the boundary $\partial M$ and the Ricci curvature satisfying ${\rm Ric}\geq(n-1)K$, where $K$ is a positive constant. If the mean curvature $H$ is positive, then for the weak solution $u$ to \eqref{1-6} with $p\geq2$, the following inequality holds:
\begin{align}\label{1-th-3-F1}
\int_{\partial M} \dfrac{1}{nH}ds\geq |M|+\frac{2(p-1)}{p}\lambda\int_{M} u^{p-1}dv.
\end{align}

\end{thm}

The corresponding version of $p=2$ with respect to inequality \eqref{1-th-3-F1} has been proved by Freitas, Roncoroni and Santos (see \cite[Theorem 1]{Freitas2024}).

\begin{thm}\label{1-th-4}(Soap Bubble Theorem)
Let $(M,g)$ be an $n$-dimensional compact Riemannian manifold with the boundary $\partial M$ and the Ricci curvature satisfying ${\rm Ric}\geq(n-1)K$, where $K$ is a positive constant. If the mean curvature $H$ is positive, then for the weak solution $u$ to \eqref{1-6} with $p\geq2$, we have
\begin{align}\label{1-th-4-F1}
\int_{\partial M} (H_0-H)|u_\nu|^{2p-2}ds\geq 0,
\end{align}
where $H_0=\frac{1}{nc}$ with $c=\frac{1}{|\partial M|}\big(|M|+\dfrac{2(p-1)}{p}\int_{M} \lambda u^{p-1}\big)$.
\end{thm}

When $p=2$ or $K=0$, as shown in Lemma \ref{2-lem-2}, then the inequality \eqref{1-th-4-F1} can be expressed as
\begin{align}\label{1-th-6-F1}
\int_{\partial M} (H_0-H)|u_\nu|^{2p-2}ds\geq0,
\end{align}
where the equality occurs if and only if $M$ is a geodesic ball and $u$ is a radial function (see \cite{Magnanini2019,Freitas2024}).

\section{Linearized $p$-operator and some lemmas}
It is well-known that the linearized operator of the $p$-Laplacian where $\nabla u\neq0$ is given (see \cite{WC2014,WCY2013,WL2016}) by
\begin{align}\label{2-1}
\mathcal{L}_{u}\eta=&{\rm div}(|\nabla u|^{p-2}A(\nabla\eta))=(|\nabla u|^{p-2}A_{ij}\eta_i)_j\notag\\
=&|\nabla u|^{p-2}\Delta \eta+(p-2)|\nabla u|^{p-4}\nabla^2\eta(\nabla u,\nabla u)\notag\\
&+(p-2)\dfrac{\left\langle\nabla u,\nabla \eta\right\rangle}{|\nabla u|^2}\Delta_p u+2(p-2)|\nabla u|^{p-4}\nabla^2u(\nabla u,\nabla \eta-\dfrac{\left\langle\nabla u,\nabla \eta\right\rangle}{|\nabla u|^2}\nabla u),
\end{align}
where $\eta\in C^2(M)$ and
$$
A=I+(p-2)\dfrac{\nabla u\otimes\nabla u}{|\nabla u|^2}.
$$
Here $I$ is the identity matrix. Let ${\mathcal{L}}_{u}^{II}$ denote the second-order part of $\mathcal{L}_u$, which is given by
\begin{align}\label{2-2}
\mathcal{L}_u^{II}\eta=&|\nabla u|^{p-2}A_{ij}\eta_{ij}\notag\\
=&|\nabla u|^{p-2}\Delta \eta+(p-2)|\nabla u|^{p-4}\nabla^2 \eta(\nabla u,\nabla u).
\end{align}
It follows that
\begin{align}\label{2-3}
\mathcal{L}_{u}\eta=&\mathcal{L}_u^{II}\eta+(p-2)\dfrac{\left\langle\nabla u,\nabla \eta\right\rangle}{|\nabla u|^2}\Delta_p u\notag\\
&+2(p-2)|\nabla u|^{p-4}\nabla^2u(\nabla u,\nabla \eta-\dfrac{\left\langle\nabla u,\nabla \eta\right\rangle}{|\nabla u|^2}\nabla u).
\end{align}
In particular, we have
$$\mathcal{L}_u(u)=(p-1)\Delta_p u$$
and
$$\mathcal{L}_u^{II} u=\Delta_p u.$$

\begin{lem}\label{2-lem-1}
Let $(M,g)$ be an $n$-dimensional compact Riemannian manifold with the Ricci curvature satisfying ${\rm Ric}\geq(n-1)K$, where $K$ is a positive constant. For the $P$-function defined in \eqref{1-8}, we have
\begin{align}\label{1-lem-1-F1}
\mathcal{L}_{u}P=&p(p-1)|\nabla u|^{2(p-2)}\left[|\nabla^{2}u|^2+{\rm Ric}(\nabla u,\nabla u)+(p-2)^2A_u^2\right]+\dfrac{p(p-1)}{n}\Delta_{p} u\notag\\
&+p(p-1)|\nabla u|^{p-2} \left\langle\nabla\Delta_{p} u,\nabla u\right\rangle+bp(p-1)^2u^{p-2}|\nabla u|^p+bp(p-1) u^{p-1}\Delta_{p} u\notag\\
&+2p(p-1)(p-2)|\nabla u|^{2(p-2)}|\nabla|\nabla u||^2
\end{align}
which holds wherever $|\nabla u|\neq0$. Here, $A_u$ is given by $A_u=\frac{\nabla^2 u(\nabla u,\nabla u)}{|\nabla u|^2}=\frac{\langle\nabla u,\nabla (|\nabla u|^2)\rangle}{2|\nabla u|^2}$.
\end{lem}

\proof
Applying \eqref{2-3}, we obtain
\begin{align}\label{2-6}
\mathcal{L}_{u}P=&\mathcal{L}_{u}^{II}P+(p-2)\frac{\langle\nabla u,\nabla P\rangle}{|\nabla u|^2}\Delta_p u+2(p-2)|\nabla u|^{p-4}\nabla^2u\Big(\nabla u,\nabla P-\frac{\langle\nabla u,\nabla P\rangle}{|\nabla u|^2}\nabla u\Big)\notag\\
=&(p-1)\mathcal{L}_{u}^{II}|\nabla u|^p+b\mathcal{L}_{u}^{II}u^p+\dfrac{p}{n}\mathcal{L}_{u}^{II}u+(p-2)\dfrac{\langle\nabla u,\nabla P\rangle}{|\nabla u|^2}\Delta_p u\notag\\
&+2(p-2)|\nabla u|^{p-4}\nabla^2u\Big(\nabla u,\nabla P-\frac{\langle\nabla u,\nabla P\rangle}{|\nabla u|^2}\nabla u\Big).
\end{align}
For the operator $\mathcal{L}_u^{II}$, we have the well-known $p$-Bochner formula (see \cite{Naber2014,Valtorta2012}): \begin{align}\label{2-5}
\frac{1}{p}\mathcal{L}_{u}^{II}(|\nabla u|^p)=&|\nabla u|^{2(p-2)}[|\nabla u|^{2-p}(\left\langle\nabla\Delta_{p} u,\nabla u\right\rangle-(p-2)A_u\Delta_{p} u)\notag\\
&+|\nabla^2 u|^2+p(p-2)A_u^2+{\rm Ric}(\nabla u,\nabla u)].
\end{align}
It is easy to see that
\begin{align}\label{2-7}
(u^p)_{ij}=p(u^{p-1}u_i)_{j}=p(p-1)u^{p-2}u_iu_j+pu^{p-1}u_{ij},
\end{align}
therefore we have
\begin{align}\label{2-9}
\mathcal{L}_{u}^{II}(u^p)=&|\nabla u|^{p-2}\Big((u^p)_{ii}+(p-2)\dfrac{(u^p)_{ij}u_iu_j}{|\nabla u|^{2}}\Big)\notag\\
=&p(p-1)^2u^{p-2}|\nabla u|^{p}+pu^{p-1}\Delta_p u.
\end{align}
On the other hand, since
\begin{align}\label{2-8}
P_{i}=&(p-1)(|\nabla u|^p)_{i}+b(u^p)_{i}+\dfrac{p}{n}u_{i}\notag\\
=&p(p-1)|\nabla u|^{p-2}u_ku_{ki}+bpu^{p-1}u_i+\dfrac{p}{n}u_i,
\end{align}
it follows that
\begin{align}\label{2-10}
\frac{\langle\nabla u,\nabla P\rangle}{|\nabla u|^2}&=\dfrac{1}{|\nabla u|^2}\Big(
p(p-1)|\nabla u|^{p-2}u_ku_{ki}u_i+bpu^{p-1}|\nabla u|^2+\dfrac{p}{n}|\nabla u|^2\Big)\notag\\
&=p(p-1)|\nabla u|^{p-2}A_u+bpu^{p-1}+\dfrac{p}{n},
\end{align}
and consequently,
\begin{align}\label{2-11}
|\nabla u|^{p-4}&\nabla^2u\Big(\nabla u,\nabla P-\frac{\langle\nabla u,\nabla P\rangle}{|\nabla u|^2}\nabla u\Big)\notag\\
&=|\nabla u|^{p-4}\Big[u_{ij}u_jP_i-\Big(p(p-1)|\nabla u|^{p-2}A_u+bpu^{p-1}+\dfrac{p}{n}\Big)u_{ij}u_iu_j\Big]\notag\\
&=p(p-1)|\nabla u|^{2(p-2)}(|\nabla|\nabla u||^2-A_u^2).
\end{align}
Inserting \eqref{2-5}, \eqref{2-9}, \eqref{2-10} and \eqref{2-11} into \eqref{2-6} yields the desired formula \eqref{1-lem-1-F1}.

\begin{lem}\label{2-lem-2}
Let $(M,g)$ be an $n$-dimensional compact Riemannian manifold with the Ricci curvature satisfying ${\rm Ric}\geq(n-1)K$, where $K$ is a positive constant.  For the solution $u$ to \eqref{1-6} with $p\geq2$, we have
\begin{align}\label{2-lem-1-F1}
\mathcal{L}_{u}P\geq(p-1)(p-2)\Big(1-\frac{1}{n}\Big)\lambda u^{p-1},
\end{align}
at points where $|\nabla u|\neq0$. In particular, when $p=2$ or $K=0$, we have $\mathcal{L}_{u}P\geq0$.
\end{lem}

\proof
Since the case where $p=2$ and $K=0$ is obviously true (see \cite{Magnanini2019,Freitas2024}), we should assume that $p>2$.

It has been proved in \cite{Ruan2024} (see the formula (12) in \cite{Ruan2024}) that
\begin{align}\label{2-13}
&|\nabla u|^{2(p-2)}\left[|\nabla^{2}u|^2+(p^2-2p+2)A_u^2\right]\notag\\
&\geq\dfrac{1}{n}(\Delta_p u)^2+\dfrac{n}{n-1}\Big(\dfrac{1}{n}\Delta_p u-(p-1)|\nabla u|^{p-2}A_u\Big)^2+2|\nabla u|^{2(p-2)}|\nabla|\nabla u||^2.
\end{align}
Substituting \eqref{2-13} into \eqref{1-lem-1-F1} yields
\begin{align}\label{2-14}
\mathcal{L}_{u}P\geq &\frac{p(p-1)}{n}(\Delta_p u)^2+\frac{np(p-1)}{n-1}\Big(\dfrac{\Delta_p u}{n}-(p-1)|\nabla u|^{p-2}A_u\Big)^2+p(p-1)|\nabla u|^{2(p-2)}{\rm Ric}(\nabla u,\nabla u)\notag\\
&+2p(p-1)^2|\nabla u|^{2(p-2)}(|\nabla|\nabla u||^2-A_u^2)+\frac{p(p-1)}{n}\Delta_{p} u
+p(p-1)|\nabla u|^{p-2}\langle\nabla\Delta_{p} u,\nabla u\rangle\notag\\
&+bp(p-1)^2u^{p-2}|\nabla u|^p+bp(p-1) u^{p-1}\Delta_{p} u\notag\\
=&\dfrac{np(p-1)}{n-1}\Big(\frac{1}{n}+\frac{\lambda u^{p-1}}{n}+(p-1)|\nabla u|^{p-2}A_u\Big)^2+2p(p-1)^2|\nabla u|^{2(p-2)}(|\nabla|\nabla u||^2-A_u^2)\notag\\
&+p(p-1)|\nabla u|^{2(p-2)}{\rm Ric}(\nabla u,\nabla u)-p(p-1)^2(\lambda-b)u^{p-2}|\nabla u|^p+p(p-1)\Big(\frac{\lambda}{n}-b\Big)u^{p-1}\notag\\
&+p(p-1)\Big(\frac{\lambda^2}{n}-b\lambda\Big)u^{2(p-1)}.
\end{align}
Applying the H\"{o}lder inequality
$$\aligned
u^{p-2}|\nabla u|^p=&\Bigg[u^{p-2}\Big(\frac{nK}{(p-1)\lambda}\Big)^{-\frac{p}{2(p-1)}}\Bigg]\Bigg[|\nabla u|^p\Big(\frac{nK}{(p-1)\lambda}\Big)^{\frac{p}{2(p-1)}}\Bigg] \\
\leq&\frac{p}{2(p-1)}\frac{nK}{(p-1)\lambda} |\nabla u|^{2(p-1)}+\frac{p-2}{2(p-1)}\Big(\frac{nK}{(p-1)\lambda}\Big)^{-\frac{p}{p-2}}u^{2(p-1)}
\endaligned$$
and using ${\rm Ric}\geq(n-1)K$, we obtain
\begin{align}\label{2-15}
&|\nabla u|^{2(p-2)}{\rm Ric}(\nabla u,\nabla u)-(p-1)(\lambda-b)u^{p-2}|\nabla u|^p+\Big(\frac{\lambda^2}{n}-b\lambda\Big)u^{2(p-1)}\notag\\
\geq&(n-1)K|\nabla u|^{2(p-1)}-(p-1)(\lambda-b)\Big[\frac{p}{2(p-1)}\frac{nK}{(p-1)\lambda} |\nabla u|^{2(p-1)}\notag\\
&+\frac{p-2}{2(p-1)}\Big(\frac{nK}{(p-1)\lambda}\Big)^{-\frac{p}{p-2}}u^{2(p-1)}\Big]+\Big(\frac{\lambda^2}{n}-b\lambda\Big)u^{2(p-1)}\notag\\
=&0
\end{align}
under the condition $\lambda>b$. Noticing again that
$$p(p-1)\Big(\frac{\lambda}{n}-b\Big)=(p-1)(p-2)\Big(1-\frac{1}{n}\Big)\lambda$$
and applying the following Cauchy Schwartz inequality:
$$A_u^2=\frac{\langle\nabla u,\nabla (|\nabla u|^2)\rangle^2}{4|\nabla u|^4}\leq |\nabla|\nabla u||^2,$$
we find that \eqref{2-14} yields
\begin{align}\label{2-16}
\mathcal{L}_{u}P\geq &\frac{np(p-1)}{n-1}\Big(\dfrac{1}{n}+\frac{\lambda}{n}u^{p-1}+(p-1)|\nabla u|^{p-2}A_u\Big)^2+(p-1)(p-2)\Big(1-\frac{1}{n}\Big)\lambda u^{p-1}\notag\\
&+2p(p-1)^2|\nabla u|^{2(p-2)}(|\nabla|\nabla u||^2-A_u^2)\notag\\
\geq&(p-1)(p-2)\Big(1-\frac{1}{n}\Big)\lambda u^{p-1}
\end{align}
and \eqref{2-lem-1-F1} follows.

If $p=2$ or $K=0$, then \eqref{2-16} becomes
\begin{align}\label{2-17}
\mathcal{L}_{u}P\geq&(p-1)(p-2)\Big(1-\frac{1}{n}\Big)\lambda u^{p-1}\notag\\
=&0,
\end{align}
which gives $\mathcal{L}_{u}P\geq0$.
We complete the proof of Lemma \ref{2-lem-2}.

\section{Proof of main results}
\subsection{Proof of Theorem \ref{1-th-1}}
The result is clear for the case of $p=2$. Therefore, we consider only the case $p>2$.

It is easy to see that on $\partial M$,
\begin{equation}\label{3-2}
A_{\nu \nu}=p-1,
\end{equation}
and from \eqref{2-8},
\begin{equation}\label{3-1}
P_\nu=p\Big((p-1)|u_\nu|^{p-2}u_{\nu\nu}+\dfrac{1}{n}\Big)u_\nu.
\end{equation}
Therefore,
\begin{align}\label{3-3}
\int_{M}\mathcal{L}_{u}Pdv=&\int_{M}{\rm div}(|\nabla u|^{p-2}A(\nabla P))dv\notag\\
=&\int_{\partial M}|u_\nu|^{p-2}A_{i\nu}P_i dv\notag\\
=&\int_{\partial M}|u_\nu|^{p-2}A_{\nu \nu}P_\nu dv\notag\\
=&p(p-1)\int_{\partial M}|u_\nu|^{p-2}\Big((p-1)|u_{\nu}|^{p-2}u_{\nu\nu}+\frac{1}{n}\Big)u_\nu ds.
\end{align}

On the boundary $\partial M$, it holds that
$$
\Delta u=(n-1)Hu_\nu+u_{\nu\nu}
$$
where $u_\nu=\langle\nabla u,\nu\rangle=-|\nabla u|$ and the equation \eqref{1-6} gives
$$
|u_\nu|^{p-2}[(p-2)u_{\nu\nu}+\Delta u]=-1-\lambda u^{p-1}.
$$
Thus,
\begin{align}\label{3-4}
|u_\nu|^{p-2}[(p-1)u_{\nu\nu}+(n-1)Hu_\nu]=-1.
\end{align}
Substituting \eqref{3-4} and
\begin{align}\label{3-5}
\int_{\partial M}|u_\nu|^{p-2}u_\nu ds&=\int_{\partial M}|\nabla u|^{p-2}u_\nu ds\notag\\
&=\int_{M}\Delta_p u dv\notag\\
&=-|M|-\lambda\int_{M}u^{p-1}dv,
\end{align}
into \eqref{3-3} yields
\begin{align}\label{3-6}
\int_{M}\mathcal{L}_{u}Pdv=&-\frac{p(p-1)(n-1)}{n}\int_{\partial M}|u_\nu|^{p-2}(1+nH|u_{\nu}|^{p-2}u_{\nu})u_\nu ds\notag\\
=&\frac{p(p-1)(n-1)}{n}\Big(|M|+\lambda \int_{M}u^{p-1}dv-\int_{\partial M}nH|u_\nu|^{2p-2}ds\Big).
\end{align}
This completes the proof.

\subsection{Proof of Theorem \ref{1-th-2}}
Note that
\begin{align}\label{3-7}
&\int_{\partial M} \dfrac{1}{nH}(1+nH|u_{\nu}|^{p-2}u_{\nu})^2ds\notag\\
&=\int_{\partial M} \dfrac{1}{nH}ds+2\int_{\partial M}|u_\nu|^{p-2}u_\nu ds+\int_{\partial M} nH|u_{\nu}|^{2p-2}ds\notag\\
&=\int_{\partial M} \dfrac{1}{nH}ds-2|M|-2\lambda\int_{M} u^{p-1}dv+\int_{\partial M} nH|u_{\nu}|^{2p-2}ds,
\end{align}
where the last equality follows from \eqref{3-5}. On the other hand, \eqref{3-6} shows
\begin{align}\label{3-8}
\int_{\partial M}nH|u_\nu|^{2p-2}ds=|M|+\lambda \int_{M}u^{p-1}dv-\frac{n}{p(p-1)(n-1)}\int_{M}\mathcal{L}_{u}Pdv.
\end{align}
Putting \eqref{3-8} into \eqref{3-7} yields
$$\aligned
&\int_{\partial M} \dfrac{1}{nH}(1+nH|u_{\nu}|^{p-2}u_{\nu})^2ds\\
=&\int_{\partial M} \dfrac{1}{nH}ds-\lambda\int_{M} u^{p-1}dv-|M|-\frac{n}{p(p-1)(n-1)}\int_{M}\mathcal{L}_{u}Pdv,
\endaligned$$
which is equivalent to the desired formula \eqref{1-th-2-F1}.

\subsection{Proof of Theorem \ref{1-th-3}}
By virtue of \eqref{1-th-2-F1}, we have
$$\aligned
\int_{\partial M}&\frac{1}{nH}ds-\lambda\int_{M} u^{p-1}dv-|M|\\
=&\frac{n}{p(p-1)(n-1)}\int_{M}\mathcal{L}_{u}Pdv+\int_{\partial M}\frac{1}{nH}(1+nH|u_{\nu}|^{p-2}u_{\nu})^2ds\\
\geq&\frac{p-2}{p}\lambda\int_{M} u^{p-1}+\int_{\partial M}\frac{1}{nH}(1+nH|u_{\nu}|^{p-2}u_{\nu})^2ds,
\endaligned$$
where the last inequality comes from \eqref{2-lem-1-F1}.
Therefore, we have
$$
\int_{\partial M}\frac{1}{nH}ds-\lambda\int_{M} u^{p-1}dv-|M|\geq\frac{p-2}{p}\lambda\int_{M} u^{p-1},
$$
and the proof of Theorem \ref{1-th-3} is finished.

\subsection{Proof of Theorem \ref{1-th-4}}
Note that
\begin{align}\label{3-9}
&\int_{\partial M}|u_\nu|^{p-2}(1+nH|u_{\nu}|^{p-2}u_{\nu})u_\nu ds\notag\\
=&\int_{\partial M}|u_\nu|^{p-2}u_\nu ds+\int_{\partial M}nH|u_\nu|^{2p-2}ds\notag\\
=&\int_{\partial M}|u_\nu|^{p-2}u_\nu ds+\int_{\partial M}n(H-H_0)|u_\nu|^{2p-2}ds+\int_{\partial M}nH_0|u_\nu|^{2p-2}ds\notag\\
=&\int_{\partial M}|u_\nu|^{p-2}u_\nu ds+\int_{\partial M}n(H-H_0)|u_\nu|^{2p-2}ds+\dfrac{1}{c}\int_{\partial M}(|u_\nu|^{p-2}u_\nu+c)^2\notag\\
&-2\int_{\partial M}|u_\nu|^{p-2}u_\nu ds-c|\partial M|\notag\\
=&-\int_{\partial M}|u_\nu|^{p-2}u_\nu ds+\int_{\partial M}n(H-H_0)|u_\nu|^{2p-2}ds+\dfrac{1}{c}\int_{\partial M}(|u_\nu|^{p-2}u_\nu+c)^2\notag\\
&-c|\partial M|.
\end{align}
We now choose the positive constant $c$ defined as
\begin{align}\label{3-10}
c&=\frac{1}{|\partial M|}\big(-\int_{\partial M} |u_\nu|^{p-2}u_\nu ds+\dfrac{p-2}{p}\int_{M} \lambda u^{p-1}dv\big)\notag\\
&=\frac{1}{|\partial M|}\Big(|M|+2(1-\dfrac{1}{p})\int_{M} \lambda u^{p-1}\Big).
\end{align}
Then \eqref{3-9} becomes
\begin{align}\label{3-11}
&\int_{\partial M}|u_\nu|^{p-2}(1+nH|u_{\nu}|^{p-2}u_{\nu})u_\nu ds\notag\\
=&-\Big(1-\frac{2}{p}\Big)\int_{M}\lambda u^{p-1} dv+\int_{\partial M}n(H-H_0)|u_\nu|^{2p-2}ds+\dfrac{1}{c}\int_{\partial M}(|u_\nu|^{p-2}u_\nu+c)^2.
\end{align}
On the other hand, \eqref{3-6} gives
\begin{align}\label{3-12}
&\int_{\partial M}|u_\nu|^{p-2}(1+nH|u_{\nu}|^{p-2}u_{\nu})u_\nu ds\notag\\
=&-\frac{n}{p(p-1)(n-1)}\int_{M}\mathcal{L}_{u}Pdv\notag\\
\leq&-\Big(1-\frac{2}{p}\Big)\lambda \int_{M}u^{p-1},
\end{align}
which implies that
\begin{align}\label{3-13}
\int_{\partial M}n(H-H_0)|u_\nu|^{2p-2}ds+\dfrac{1}{c}\int_{\partial M}(|u_\nu|^{p-2}u_\nu+c)^2\leq 0.
\end{align}
Then the inequality \eqref{1-th-4-F1} follows from \eqref{3-11} and \eqref{3-13}.

\bibliographystyle{Plain}

\end{document}